\documentclass{amsart} 
\input diagrams
\numberwithin{equation}{section}
\newtheorem{theorem}[equation]{Theorem}
\newtheorem{proposition}[equation]{Proposition}
\newtheorem{corollary}[equation]{Corollary}

\newtheorem{lemma}[equation]{Lemma}

\theoremstyle{remark}
\newtheorem{definition}[equation]{Definition}
\newtheorem{example}[equation]{Example}
\newtheorem*{remark}{Remark}

\newarrow{Mapsto}{|}{-}{-}{-}{>}
\newarrow{Equal}{=}{=}{=}{=}{=}
\DeclareMathOperator{\rk}{rk}
\DeclareMathOperator{\ad}{ad}
\DeclareMathOperator{\End}{End}

\begin{document}
\title{Stable bundles on non-K\"ahler elliptic surfaces}

\author{Vasile Br\^{\i}nz\u{a}nescu}
\address{Institute of Mathematics "Simion Stoilow",
Romanian Academy, P.O.Box 1-764, RO-70700,
Bucharest, Romania}
\email{Vasile.Brinzanescu@imar.ro} 

\thanks{The first author was partially supported by Swiss NSF contract SCOPES
2000-2003, No.7 IP 62615 and by contract CERES 39/2002-2004}

\author{Ruxandra Moraru}
\address{Department of Mathematics and Statistics, Burnside Hall,
McGill University, 805 Sherbrooke Street West,
Montreal, Quebec, Canada, H3A 2K6}
\email{moraru@math.mcgill.ca}


\begin{abstract}
In this paper, we study the moduli spaces $\mathcal{M}_{\delta,c_2}$ of stable rank-2 
vector bundles on non-K\" ahler elliptic surfaces, thus giving a classification 
these bundles; 
in the case of Hopf and Kodaira surfaces,
these moduli spaces admit the structure of an algebraically 
completely integrable Hamiltonian system.
\end{abstract}

\maketitle

\section{Introduction}
\indent 
Vector bundles on elliptic fibrations have been extensively studied over the past fifteen years;
in fact, there is by now a well-understood theory for projective elliptic surfaces 
(see, for example, \cite{D,F1,FMW}).
However, not very much is known about the non-K\"{a}hler case.
In this article, we partly remedy this problem by examining the stability properties 
of holomorphic rank-2 vector bundles 
on non-K\" ahler elliptic surfaces; their existence and classification 
are investigated in \cite{Brinzanescu-Moraru1,Brinzanescu-Moraru2}. 
One of the motivations for the study of bundles on non-K\"{a}hler elliptic fibrations 
comes from recent developments in superstring theory,
where six-dimensional non-K\"{a}hler manifolds occur in the context of type IIA string
compactifications with non-vanishing background $H$-field;
in particular, all the non-K\"{a}hler examples appearing in the physics litterature so far
are non-K\"{a}hler principal elliptic fibrations 
(see \cite{BBDG,CCDLMZ,GP} and the references therein).
The techniques developed here and in 
\cite{Brinzanescu-Moraru1,Brinzanescu-Moraru2} can also be used to
study holomorphic vector bundles of arbitrary rank on higher dimensional 
non-K\"{a}hler elliptic and torus fibrations. 

A minimal non-K\"{a}hler elliptic surface $X$ is a Hopf-like surface that
admits a holomorphic fibration $\pi: X \rightarrow B$, over a smooth connected compact curve $B$, 
whose smooth fibres are isomorphic to a fixed smooth elliptic curve $T$; 
the fibration $\pi$ can have at most a finite number of singular fibres, 
in which case, they are isogeneous to multiples of $T$.
More precisely, if the surface $X$ does not have multiple fibres, 
then it is the quotient of a complex surface
by an infinite cyclic group (see example \ref{degree of line bundles}); 
multiple fibres can then be introduced by performing a finite number
of logarithmic transformations on its relative Jacobian $J(X)$.
To study bundles on $X$, a natural operation is restriction to the smooth fibres of $\pi$; 
this gives rise to an important invariant, 
called the {\em spectral curve} or {\em cover}, 
which is an effective divisor on $J(X)$ that encodes the holomorphic type of 
the bundle over each smooth fibre of $\pi$. 

Consider the moduli space ${\mathcal M}_{\delta,c_2}$ of stable holomorphic 
rank-$2$ vector bundles on $X$ with fixed determinant $\delta$ and second Chern class $c_2$.
Note that the line bundle $\delta$ on $X$ induces an involution $i_\delta$ of the 
relative Jacobian $J(X)$;
furthermore, the spectral cover of any bundle in ${\mathcal M}_{\delta,c_2}$ is invariant
with respect to this involution, 
thus descending to an effective divisor on the ruled surface $\mathbb{F}_\delta := J(X)/i_\delta$
called the {\em graph} of the bundle.
We can then define a map 
\[ G: {\mathcal M}_{\delta,c_2} \rightarrow \rm{Div}(\mathbb{F}_\delta) \]
that associates to each stable vector bundle its graph in $\rm{Div}(\mathbb{F}_\delta)$,
called the {\em graph map}. In \cite{B-H,Moraru}, the stability properties of vector bundles
on Hopf surfaces were studied by analysing the image and the fibres of this map;
in particular, it was shown \cite{Moraru} that the moduli spaces admit a natural Poisson
structure with respect to which the graph map is a Lagrangian fibration
whose generic fibre is an abelian variety:
the map $G$ admits an algebraically completely integrable system structure.
In this paper, we adopt this approach to study stable vector bundles on arbitrary 
non-K\"{a}hler surfaces. 
 
The article is organised as follows. 
We begin with a brief review of some existence and classification results for
holomorphic vector bundles on non-K\"{a}hler elliptic surfaces that were proven in
\cite{Brinzanescu-Moraru1,Brinzanescu-Moraru2}.
In the third section, we obtain explicit conditions for the stability of rank-2 vector bundles:
we show that unfiltrable bundles are always stable
and then classify the destabilising bundles of filtrable bundles.
The moduli spaces ${\mathcal M}_{\delta,c_2}$ are studied in the last section. 
We first prove that these spaces are smooth on an open dense subset consisting
of vector bundles that are regular on every fibre of $\pi$
(on a smooth elliptic curve, a bundle of degree zero is said to be {\em regular}
if its group of automorphisms is of the smallest possible dimension).
However, for Hopf and Kodaira surfaces, the moduli spaces are also smooth at points 
that are not regular; in this case, the moduli are smooth complex manifolds of dimension
$4c_2 - c_1^2(\delta)$.
Then, we determine the image of the graph map; for simplicity,
we focus our presentation on non-K\"{a}hler elliptic surfaces without multiple fibres,
but similar results hold in the multiple fibre case.
Furthermore, we give an explicit description of the fibres of the graph map,
which follows immediately from the classification results of \cite{Brinzanescu-Moraru2,Moraru}
and the stability conditions of the third section; 
in particular, the generic fibre at a graph $\mathcal{G} \in Div(\mathbb{F}_\delta)$
is isomorphic to a finite number of copies of a Prym variety
associated to $\mathcal{G}$. 
We conclude by noting that for Kodaira surfaces the graph map 
is also an algebraically completely integrable Hamiltonian system, with respect
to a given symplectic structure on ${\mathcal M}_{\delta,c_2}$. 
\vspace{.05in}

\noindent
{\bf Acknowledgements.} 
The first author would like to express his gratitude to the 
Max Planck Institute of Mathematics for its hospitality and 
stimulating atmosphere; part of this paper was prepared during his stay at
the Institute. 
The second author would like to thank Jacques Hurtubise for his generous 
encouragement and support during the completion of this paper;
she would also like to thank Ron Donagi and Tony Pantev for valuable discussions, and the
Department of Mathematics at the University of Pennsylvania 
for their hospitality,
during the preparation of part of this article.

\section{Holomorphic vector bundles}
\label{vector bundles}

Let $X\stackrel{\pi}{\rightarrow}B$ be a minimal non-K\"ahler elliptic 
surface, with $B$ a smooth compact connected curve; it is well-known that 
$X \stackrel{\pi}{\rightarrow} B$ is a quasi-bundle over $B$, that is, all the smooth fibres 
are pairwise isomorphic and the singular fibres are multiples of elliptic 
curves \cite{Kod,Br}.
Let $T$ be the general fibre of $\pi$, which is an elliptic curve, and
denote its dual $T^*$ (we fix a non-canonical identification 
$T^*:= \mbox{Pic}^0(T)$); in this case, the relative Jacobian of 
$X\stackrel{\pi}{\rightarrow}B$  is simply 
\[ J(X)=B\times T^*\stackrel{p_1}{\rightarrow}B \]
(see, for example, \cite{Kod,BPV,Br}) and 
$X$ is obtained from $J(X)$ by a finite number 
of logarithmic transformations \cite{Kod,BPV,BrU}.
In addition, if the fibration $\pi$ has multiple fibres, then one can associate to 
$X$ a principal $T$-bundle $\pi': X' \rightarrow B'$ over an $m$-cyclic covering 
$\varepsilon: B' \rightarrow B$, where the integer $m$ depends on the multiplicities
of the singular fibres;
note that the map $\varepsilon$ induces natural $m$-cyclic coverings
$J(X') \rightarrow J(X)$ and $\psi :X'\rightarrow X$. 

To study bundles on $X$, one of our main tools is 
restriction to the smooth fibres of the fibration $\pi: X \rightarrow B$. 
It is important to point out that since $X$ is non-K\"{a}hler,
the restriction of {\em any} vector bundle on $X$ to a smooth fibre of $\pi$ {\em always} 
has trivial first Chern class \cite{Brinzanescu-Moraru1}.
Therefore, a vector bundle $E$ on $X$ is semistable on the generic
fibre of $\pi$; in fact, its restriction to a fibre $\pi^{-1}(b)$
is unstable on at most an isolated set of points $b \in B$; 
these isolated points are called the {\em jumps} of the bundle.
Furthermore, there exists a divisor $S_E$ in the relative Jacobian of $X$, 
called the {\em spectral curve} or {\em cover} of the bundle, 
that encodes the isomorphism
class of the bundle $E$ over each smooth fibre of $\pi$;
for a detailed description of this divisor, we refer the reader to 
\cite{Brinzanescu-Moraru1}.
We should note that, if the fibration $\pi$ has multiple fibres, then the spectral 
cover $S_E$ of $E$ is actually defined as the projection in $J(X)$ of the spectral cover 
$S_{\psi^\ast E} \subset J(X')$ of $\psi^\ast E$, where $\psi: X' \rightarrow X$ is 
the $m$-cyclic covering defined above.

\subsection{Line bundles} 
\label{line bundles}

The spectral cover of a line bundle $L$
on $X$ is a section $\Sigma$ of $J(X)$ such that the restriction of $L$ to
any smooth fibre $\pi^{-1}(b)$ of $\pi$ is isomorphic to the line bundle $\Sigma_b$
of degree zero on $T \cong \pi^{-1}(b)$.
Conversely, given any section $\Sigma$ of $J(X)$, there exists at least one line
bundle on $X$ with spectral cover $\Sigma$ \cite{Brinzanescu-Moraru1}.
Before giving a classification of line bundles on $X$, we fix some notation.
Suppose that $\pi$ has a multiple fibre $mF$ over the point $b$ in $B$; the line bundle 
associated to the divisor $F$ of $X$ is then such that 
$(\mathcal{O}_X(F))^m = \mathcal{O}_X(mF) = \pi^\ast\mathcal{O}_B(b)$. 
Let $P_2$ be the subgroup of $\text{Pic}(X)$ generated by $\pi^\ast \text{Pic}(B)$ and
the $\mathcal{O}_X(T_i)'$, where $m_1T_1, \dots, m_rT_r$ are the multiple fibres (if any) of $X$;
we then have \cite{Brinzanescu-Moraru1}:  

\begin{proposition}
Let $\Sigma$ be a section of $J(X)$. Then, the set of all line bundles on $X$ 
with spectral cover $\Sigma$ is a principal homogeneous space over $P_2$.
\end{proposition} 

\subsection{Rank-2 vector bundles}
\label{rank-2 vector bundles}

Consider a rank-2 vector bundle $E$ on $X$; its {\em discriminant} is then defined as
\[ \Delta(E) := \frac{1}{2} \left( c_2(E) - \frac{c_1(E)^2}{4} \right). \]
In this case, the spectral curve of $E$ is a divisor $S_E$ in $J(X)$ of the form
\[ S_E := \left( \sum_{i=1}^k \{ x_i \} \times T^\ast \right) + \overline{C}, \]
where $\overline{C}$ is a bisection of $J(X)$ and $x_1, \cdots, x_k$ are points in $B$
that correspond to the jumps of $E$. 
Let $\delta$ be the determinant line bundle of $E$. 
It then defines the following involution on the relative Jacobian $J(X) = B \times T^\ast$ of $X$:
\[ \begin{array}{rcl}
i_\delta : J(X) & \rightarrow & J(X)\\
(b,\lambda) & \mapsto & (b,\delta_b \otimes \lambda^{-1}),
\end{array} \] 
where $\delta_b$ denotes the restriction of $\delta$ to the fibre 
$T_b = \pi^{-1}(b)$.
By construction, the spectral curve $S_E$ of $E$ is invariant with respect to
this involution; in particular, the pair of points lying on the bisection 
$\overline{C}$ over $b$ is of the form $\{ \lambda_b, \delta_b \otimes \lambda_b^{-1} \}$,
where $\lambda_b$ and $\delta_b \otimes \lambda_b^{-1}$ are 
the subline bundles of $E|_{\pi^{-1}(b)}$.
Finally, note that the quotient of $J(X)$ by the involution $i_\delta$ is a ruled surface 
$\mathbb{F}_\delta := J(X)/i_\delta$ over $B$;
let $\eta: J(X) \rightarrow \mathbb{F}_\delta$ be the canonical map.
The spectral cover $S_E$ of $E$ then descends to 
divisor on $\mathbb{F}_\delta$
of the form 
\[ \mathcal{G}_E := \sum_{i=1}^k f_i + A,\]
where $f_i$ is the fibre of the ruled surface $\mathbb{F}_\delta$ over 
the point $x_i$ and
$A$ is a section of the ruling such that $\eta^\ast A = \overline{C}$. 

\subsubsection{Bundles without jumps}

We begin with some properties of filtrable bundles without jumps.
Let $E$ be a rank-2 vector bundle on $X$ with determinant $\delta$,
and spectral cover $\left( \Sigma_1 + \Sigma_2 \right)$, 
where $\Sigma_1$ and $\Sigma_2$ are
sections of $J(X)$;
there exists a line bundle $\mathcal{D}$ on $X$ associated to $\Sigma_1$
such that $E$ is given by an extension 
\begin{equation}\label{line bundle extension} 
0 \rightarrow \mathcal{D} \rightarrow E \rightarrow 
{\mathcal{D}}^{-1} \otimes \delta \rightarrow 0. 
\end{equation}
Consequently, 
\begin{equation}\label{Delta(E)}
\Delta(E) = - \frac{1}{8} \left( c_1(\delta) - 2c_1(\mathcal{D}) \right)^2.
\end{equation}
Given the above considerations, we have the following results.
\begin{lemma}\label{multiple section}
If $\Sigma_1 = \Sigma_2$, then $\Delta(E) = 0$. 
Furthermore, the extension \eqref{line bundle extension}
either splits on every fibre of $\pi$ or else it splits on at most a finite number 
of fibres. 
\end{lemma}
\begin{proof} 
Note that $c_1(\mathcal{D}) = c_1(\mathcal{D}^{-1} \otimes \delta)$
because $\Sigma_1 = \Sigma_2$;
referring to \eqref{Delta(E)}, we then have $\Delta(E) = 0$.
Suppose that there exists at least one fibre $T_{b_0}$ of $\pi$ over which the extension is 
non-trivial; therefore, $h^1(T_{b_0}, \mathcal{D}^{-1} \otimes E) = 1$. But if the extension splits
over the fibre $T_b$, then $h^1(T_b, \mathcal{D}^{-1} \otimes E) = 2$.
The upper semi-continuity of the map $b \mapsto h^1(T_b, \mathcal{D}^{-1} \otimes E)$
thus implies that $h^1(T_b, \mathcal{D}^{-1} \otimes E) = 1$ for generic $b$. 
\end{proof}
\begin{lemma}\label{not multiple section}
If $\Sigma_1 \neq \Sigma_2$, then $|\Sigma_1 \cap \Sigma_2| = 4\Delta(E)$.
In addition, the extension \eqref{line bundle extension}
splits globally whenever $\Delta(E) = 0$.
\end{lemma}
\begin{proof}
Since $\Sigma_1 \neq \Sigma_2$, the sheaf 
$\pi_\ast(\mathcal{D}^{-2} \otimes \delta)$ vanishes and the first direct image sheaf
$R^1\pi_\ast(\mathcal{D}^{-2} \otimes \delta)$ is a skyscraper sheaf supported
on the points of $\Sigma_1 \cap \Sigma_2$. Therefore, 
$c_1(R^1\pi_\ast(\mathcal{D}^{-2} \otimes \delta)) = |\Sigma_1 \cap \Sigma_2|$ and by
Grothendieck-Riemann-Roch,
\[ |\Sigma_1 \cap \Sigma_2| = - \frac{1}{2} \left( c_1(\delta) - 2c_1(\mathcal{D}) \right)^2, \]
which is equal to $4\Delta(E)$ by \eqref{Delta(E)}.
Consequently, if $\Delta(E) = 0$, then $\Sigma_1 \cap \Sigma_2  = \emptyset$; 
in this case, the extension \eqref{line bundle extension} splits on every fibre of $\pi$
and $R^1\pi_\ast({\mathcal{D}}^2 \otimes \delta^{-1}) = 0$.
Hence, the Leray spectral sequence gives
$H^1(X,{\mathcal{D}}^2 \otimes \delta^{-1}) = 0$ and the extension
splits globally.
\end{proof}

We have seen that we can associate to any rank-2 vector bundle on $X$ a bisection in $J(X)$.
Conversely, given any bisection of $J(X)$, there exists at least one rank-2 vector 
bundle on $X$ associated to it; if the bisection is smooth, 
the bundles that correspond to it are classified as follows  
(see \cite{Brinzanescu-Moraru2} for precise statements).
\begin{theorem}\label{classification of regular bundles}
Fix a line bundle $\delta$ on $X$ and its associated involution $i_\delta$ of $J(X)$. 
Let $C$ be a smooth bisection of $J(X)$ that
is invariant with respect to this involution; it is then a double
cover of $B$ of genus $4\Delta(2,c_1,c_2) + 2g - 1$.
The set of all rank-2 vector bundles on $X$ with spectral cover $C$ and 
determinant $\delta$ is then parametrised by a finite number of copies of the 
Prym variety $Prym(C/B)$ associated to the double cover $C \rightarrow B$.
In particular, if $X$ does not have multiple fibres, then this set
of bundles consists of a single copy of $Prym(C/B)$. 
\end{theorem}

\subsubsection{Bundles with jumps}

Consider a rank-2 vector bundle $E$ on $X$ 
with determinant $\delta$ that has a jump of multiplicity $\mu$ over the smooth fibre 
$T = \pi^{-1}(x_0)$.
The restriction of $E$ to the fibre $T$ is then of the 
form $\lambda \oplus (\lambda^\ast \otimes \delta_{x_0})$, for some 
$\lambda \in \text{Pic}^{-h}(T)$, $h>0$;
the integer $h$ is called the {\it height} of the jump at $T$.
Moreover, up to a multiple of the identity, 
there is a {\em unique} surjection $E|_T \rightarrow \lambda$,
which defines a canonical elementary modification of $E$ that we denote $\bar{E}$; 
this elementary modification is called {\em allowable} \cite{F2}.
Therefore, we can associate to $E$ a finite sequence
$\{ \bar{E}_1, \bar{E}_2, \dots , \bar{E}_l \}$
of allowable elementary modifications such that $\bar{E}_l$ is
the only element of the sequence that does not have a jump at $T$.

Let us now assume that $\pi$ has a multiple fibre $m_0T_0$. 
One can then associate to $X$ an elliptic quasi-bundle
$\pi': X' \rightarrow B'$, over an $m_0$-cyclic covering $\varepsilon: B' \rightarrow B$,
such that $T'_0 := \psi^{-1}(T_0) \subset X'$
is a smooth fibre of $\pi'$, where $\psi: X' \rightarrow X$ is the
$m_0$-cyclic covering induced by $\varepsilon$.  
Given this, we say that $E$ has a jump over $T_0$
if and only if the restriction of $\psi^\ast E$ to the fibre $T'_0$ is unstable. 
Naturally, the height and multiplicity of the jump of $E$ over $T_0$
are defined as the height and multiplicity of the jump of $\psi^\ast E$ over $T'_0$.
We can now define the following important invariants.

\begin{definition}\label{jumping sequence - definition}
Let $T$ be a smooth fibre of $\pi$. Suppose that the vector bundle $E$ 
has a jump over $T$ and consider the corresponding sequence of allowable elementary modifications
$\{ \bar{E}_1, \bar{E}_2, \dots , \bar{E}_l \}$. 
The integer $l$ is called
the {\it length} of the jump at $T$. The {\it jumping sequence} of $T$ is defined
as the set of integers $\{ h_0, h_1, \dots , h_{l-1} \}$, where $h_0 = h$
is the height of $E$ and $h_i$ is the height of
$\bar{E}_i$, for $0 < i \leq l-1$.

If the vector bundle $E$ has a jump over a multiple fibre $m_0T_0$ of $\pi$, we define
the length and jumping sequence of $T_0$ to be the length and jumping sequence of 
the jump of $\psi^\ast E$ over the smooth fibre $T'_0 = \psi^{-1}(T_0)$ of $\psi$, 
where $\psi: X' \rightarrow X$ is the $m_0$-cyclic covering defined above.
\end{definition}

Note that if a vector bundle $E$ jumps over a smooth fibre 
$T$ of $\pi$, with  multiplicity $\mu$ and jumping sequence $\{h_0, \dots, h_{l-1}\}$,
then $\mu = \sum_{i=1}^{l-1} h_i$.
For a detailed description of jumps, we refer the reader to 
\cite{Moraru,Brinzanescu-Moraru2}; moreover, the basic properties of elementary
modifications can be found, for example, in \cite{F2}.
We finish this section by stating the following existence result \cite{Brinzanescu-Moraru2}: 
if $X$ does not have multiple fibres and $S$ is a spectral cover in $J(X)$ 
(that may have vertical components), then one can associate to $S$ at least one 
rank-2 vector bundle on $X$.

\section{Stable rank two bundles}
\subsection{Degree and stability}
\label{degree}
\indent
The degree of a vector bundle can be defined on any compact
complex manifold $M$. Let $d = \dim_\mathbb{C} M$. 
A theorem of Gauduchon's \cite{Gauduchon} states that any hermitian metric on $M$ 
is conformally equivalent to a metric, 
called a {\it Gauduchon metric}, whose associated (1,1) form 
$\omega$ satisfies $\partial \bar{\partial} \omega^{d-1} = 0$.
Suppose that $M$ is endowed with such a metric and let $L$ be a holomorphic 
line bundle on $M$. The {\it degree of $L$ with respect to $\omega$} is 
defined \cite{Buchdahl}, up to a constant factor, by
\[ \deg{L} := \int_M F \wedge \omega^{d-1},\]
where $F$ is the curvature of a hermitian connection on $L$, compatible with 
$\bar{\partial}_L$. Any two such forms $F$ differ by a 
$\partial \bar{\partial}$-exact form. Since 
$\partial \bar{\partial} \omega^{d-1} = 0$, 
the degree is independent of the choice of connection and is therefore
well defined. This notion of degree is an extension of the K\"{a}hler case. 
If $M$ is K\"{a}hler, we get the usual topological degree defined
on K\"{a}hler manifolds; but in general, this degree is 
not a topological invariant, for it 
can take values in a continuum (see below).

Having defined the degree of holomorphic line bundles, we define the 
{\it degree} of a torsion-free coherent sheaf $\mathcal{E}$ on $M$ by
\[ \deg(\mathcal{E}) := \deg(\det{\mathcal{E}}), \]
where $\det{\mathcal{E}}$ is the determinant line bundle of $\mathcal{E}$,
and the {\it slope of $\mathcal{E}$} by 
\[ \mu(\mathcal{E}) := \deg(\mathcal{E})/\text{rk}(\mathcal{E}).\]
The notion of stability then exists for any compact complex manifold:
\vspace{.1in}
\newline
{\em A torsion-free coherent sheaf $\mathcal{E}$ on $M$ is {\em stable} 
if and only if for every coherent subsheaf $\mathcal{S} \subset \mathcal{E}$ 
with $0 < \text{rk}(\mathcal{S}) < \text{rk}(\mathcal{E})$, we have
$\mu(\mathcal{S}) < \mu(\mathcal{E})$.}
\begin{remark} With this definition of stability, many of the 
properties from the K\"{a}hler case hold. In particular, 
all line bundles are stable;
for rank two vector bundles on a surface, 
it is sufficient to verify stability with respect to line bundles.
Finally, if a vector bundle $E$ is stable, then $H^0(M,End(E)) = \mathbb{C}$. 
\end{remark}

\begin{example}\label{degree of line bundles}
Let $X \stackrel{\pi}{\rightarrow} B$  
be a non-K\"{a}hler principal elliptic bundle over a curve $B$ of genus $g$ and with fibre $T$.
The surface $X$ is then isomorphic to a quotient of the form
\[ X = \Theta^\ast / \langle \tau \rangle ,\]
where $\Theta$ is a line bundle on $B$ with positive Chern class $d$,
$\Theta^\ast$ is the complement of the zero section in the total space of 
$\Theta$, and $\langle \tau \rangle$ is the multiplicative cyclic group 
generated by a fixed complex number $\tau \in \mathbb{C}$, with 
$| \tau | > 1$. 
In this case, the degree of torsion line bundles can be 
computed explicitly (for details, see \cite{T}).
Every line bundle $L \in \text{Pic}^{\tau}(X)$ decomposes uniquely 
as $L = H \otimes L_\alpha$, for $H \in \cup_{i=0}^{d-1} \text{Pic}(B)$
and $\alpha \in \mathbb{C}^\ast$. Taking into account this
decomposition, the degree of $L$ is given by
\[ \deg{L} = c_1(H) - \frac{d}{\ln|\tau|}\ln|\alpha|. \]
In particular, $\deg(\pi^\ast H)= \deg{H}$, for all $H \in \text{Pic}(B)$. 

We end this example by observing that if $X$ has a multiple fibre $m_0T_0$, then we have
$\deg(\mathcal{O}_X(T_0)) = 1/m_0$.
\end{example}

\subsection{Stable vector bundles}
\label{stable bundles}

Let $X \stackrel{\pi}{\rightarrow} B$ be a non-K\"{a}hler elliptic surface  
with multiple fibres $m_1T_1$, \dots , $m_rT_r$ (if any);
the canonical bundle of $X$ is then
$K_X = \pi^\ast K_B \otimes \mathcal{O}_X \left( \sum_{i=1}^r (m_i - 1) T_i \right)$,
giving
$\omega_{X/B} = \mathcal{O}_X \left( \sum_{i=1}^r (m_i - 1) T_i \right)$
as the dualising sheaf of $\pi$.
Note that $\deg{\omega_{X/B}}  = r - \sum_{i=1}^r 1/m_i \geq 0$ 
(see  example \ref{degree of line bundles}).
Fix a rank-2 vector bundle $E$ on $X$ and let $\delta$
be its determinant line bundle; 
there exists a sufficient condition on the spectral cover of $E$ 
that ensures its stability:

\begin{proposition}\label{spectral cover of stable bundles}
Suppose that the spectral cover of $E$ includes an irreducible bisection 
$\overline{C}$ of $J(X)$. Then $E$ is stable.
\end{proposition}
\begin{proof}
Suppose that there exists a line bundle $\mathcal{D}$ on $X$ that maps into $E$.
After possibly tensoring $\mathcal{D}$ by the pullback of a suitable line 
bundle on $B$, the rank-2 bundle $E$ is then given as an extension
\begin{equation}\label{extension}
0 \rightarrow \mathcal{D} \rightarrow E \rightarrow 
\mathcal{D}^{-1} \otimes \delta \otimes I_Z \rightarrow 0,
\end{equation}
where $Z \subset X$ is a locally complete intersection of codimension 2. 
In fact, $Z$ is the set of points $\{ x_1, \dots , x_k \}$ corresponding 
to the fibres $\pi^{-1}(x_i)$ over which $E$ is unstable.
Let $\Sigma_1$ and $\Sigma_2$ be the sections of $J(X)$ determined by the line
bundles $\mathcal{D}$ and $\mathcal{D}^{-1} \otimes \delta$, respectively.
The extension \eqref{extension} then implies 
$\overline{C} = \Sigma_1 + \Sigma_2$. 
\end{proof}

Consequently, the spectral covers of unstable bundles 
include bisections of the form 
$\overline{C} = \left( \Sigma_1 + \Sigma_2 \right)$, where $\Sigma_1$ and $\Sigma_2$ are
sections of the Jacobian surface.

\begin{proposition}\label{destab. bundles}
Suppose that the spectral cover of $E$ is given by
\[ \left( \sum_{i=1}^k \{ x_i \} \times T^\ast \right) + 
\left( \Sigma_1 + \Sigma_2 \right). \] 
Then, there exist line bundles $K_1$ and $K_2$ on $X$  
(corresponding to the sections $\Sigma_1$ and $\Sigma_2$, respectively) 
such that the set 
of all line bundles that map non-trivially to $E$ is given by
\[ \left\{ K_j \otimes \pi^\ast H \otimes 
\mathcal{O}_X \left(\sum_{i=1}^r a_i T_i \right) \ : \ 
\mbox{$H \in \text{\em Pic}^{\leq 0}(B)$ and $a_i \leq 0$} \right\}.\] 
Also, $E$ is stable if and only if $\deg{K_1}$ and $\deg{K_2}$ are both
smaller than $\deg{\delta}/2$.
Note that if $\Sigma_1 = \Sigma_2$, then $K_1 = K_2$.
 
The line bundles $K_1$ and $K_2$ are called the {\em destabilising line bundles of $E$}.
\end{proposition}
\begin{proof}
Let $\mathcal{D}$ be a line bundle that corresponds to the section $\Sigma_1$
and suppose that there exists a non-trivial map $\mathcal{D} \rightarrow E$.
We begin by assuming that $E$ is regular on the generic fibre of $\pi$.
In this case, the direct image sheaf $\pi_\ast(\mathcal{D}^{-1} \otimes E)$
is a line bundle on $B$, say $L$, of positive degree. 
Set $K = \mathcal{D} \otimes (\pi^\ast L)^{-1}$;
then, $K$ restricts to ${\Sigma_1}_b$ over the smooth fibres $\pi^{-1}(b)$ of $B$
and $\pi_\ast({K}^{-1} \otimes E) \cong \mathcal{O}_B$.
However, any line bundle $\mathcal{D}'$ on $X$ corresponding to $\Sigma_1$
can be written as 
$K \otimes \pi^\ast H \otimes \mathcal{O}_X \left( \sum_{i=1}^r a_i T_i \right)$, 
for some $H \in \text{Pic}(B)$ and integers $0 \leq a_i \leq m_i - 1$.
Moreover, one can easily show that
$\pi_\ast \left(\mathcal{F} \otimes \mathcal{O}_X \left( \sum_{i=1}^r a_i T_i \right) \right) 
= \pi_\ast(\mathcal{F})$,  
for any locally free sheaf $\mathcal{F}$ on $X$.
Consequently, if $\mathcal{D}'$ also maps into $E$, then the line bundle
$\pi_\ast({\mathcal{D}'}^{-1} \otimes E) \cong H^{-1}$ has a non-trivial section,
implying that $H \in \text{Pic}^{\leq 0}(B)$.
Note that $\pi_\ast({K'}^{-1} \otimes E) \cong \mathcal{O}_B$
for any line bundle $K'$ of the form 
$K \otimes \mathcal{O}_X \left( \sum_{i=1}^r b_i T_i \right)$,
where $0 \leq b_i \leq m_i - 1$. 
But, the destabilising bundle $K_1$ is the line bundle associated to $\Sigma_1$
of maximal degree that maps into $E$; we therefore set 
$K_1 =  K \otimes \mathcal{O}_X \left( \sum_{i=1}^r (m_i - 1) T_i \right) = K \otimes \omega_{X/B}$.
Clearly, any line bundle corresponding to $\Sigma_1$ that maps into $E$ is 
can be written as 
$K_1 \otimes \pi^\ast H \otimes \mathcal{O}_X \left( \sum_{i=1}^r a_i T_i \right)$, 
for $H \in \text{Pic}^{\leq 0}(B)$ and integers $a_i \leq 0$. 

We now assume that $E$ is not regular on the generic fibre of $\pi$.
The direct image sheaf $\pi_\ast(\mathcal{D}^{-1} \otimes E)$
is then a rank-2 vector bundle on $B$, say $\mathcal{F}$;
it must have a subline bundle $L$ such that $\mathcal{F}/L$ is torsion free.
If we set 
$K_1 = \mathcal{D} \otimes (\pi^\ast L^{-1})^{-1}
\otimes \omega_{X/B}$, 
then $\pi_\ast({K_1}^{-1} \otimes E)$ has a nowhere vanishing section and,
as above, any line bundle induced by $\Sigma_1$ that maps into $E$ is 
of the required form.
\end{proof}

In fact, the destabilising line bundles of filtrable bundles without jumps can 
be described explicitely as follows:

\begin{proposition}\label{max. destab. bundles}
Let $E$ be a holomorphic rank-2 vector bundle on $X$ with invariants $\det(E) = \delta$,
$c_2(E) = c_2$, and spectral cover $\left( \Sigma_1 + \Sigma_2 \right)$, 
where $\Sigma_1$ and $\Sigma_2$ are
sections of $J(X)$.
Let $K_1$ be the destabilising line bundle of $E$ induced by $\Sigma_1$;
there is an extension
\begin{equation}\label{line bundle extension2}
0 \rightarrow K_1 \rightarrow E \rightarrow 
{K_1}^{-1} \otimes \delta \rightarrow 0.
\end{equation}

(i) If the extension is trivial on every fibre of $\pi$,
then the second destabilising bundle of $E$ is $K_2 = {K_1}^{-1} \otimes \delta \otimes \omega_{X/B}$. 

(ii) Suppose that $\Sigma_1 = \Sigma_2$ and that
the extension splits on a finite number $n \geq 0$ of fibres of $\pi$.
In this case, ${K_1}^2 = \delta \otimes \pi^\ast (H_+) 
\otimes \omega_{X/B} $, 
where $H_+$ is a line bundle of degree $n$ on $B$ that is trivial whenever $n=0$.

(iii) If $\Sigma_1 \neq \Sigma_2$ and the extension is non-trivial on a finite number 
$n \leq 4\Delta(E)$ 
of fibres, then the second destabilising bundle of $E$ is 
$K_2 = {K_1}^{-1} \otimes  \delta \otimes \pi^\ast (H_-) \otimes \omega_{X/B}$, 
for some line bundle $H_-$ on $B$ of negative degree $-n$.
\end{proposition}

\begin{proof}
Let us first assume that the extension 
\eqref{line bundle extension2} splits on every fibre of $\pi$; 
the line bundle $\pi_\ast(K_1 \otimes \delta^{-1} \otimes E)$
has a nowhere vanishing global section and is thus trivial.
The second destabilising bundle
of $E$ is then $K_2 = {K_1}^{-1} \otimes \delta  \otimes \omega_{X/B}$, proving (i). 

Next, we suppose that $\Sigma_1 = \Sigma_2$. The restriction of
${K_1}^2 \otimes \delta^{-1}$ is therefore trivial every fibre of $\pi$, implying that 
${K_1}^2 \otimes \delta^{-1}= \pi^\ast (H_+) 
\otimes \mathcal{O}_X \left( \sum_{i=1}^r a_i T_i \right)$ for a line bundle $H_+$ on $B$
and integers $0 \leq a_i \leq m_i - 1$;
tensoring the exact sequence \eqref{line bundle extension2}
by ${K_1}^{-1}$ and pushing down to $B$, we obtain a new exact sequence:
\begin{equation}\label{long exact} 
0 \rightarrow {H_+}^{-1} \rightarrow \mathcal{O}_B \rightarrow 
R^1\pi_\ast({K_1}^{-1} \otimes E) \rightarrow{H_+}^{-1} \rightarrow 0. 
\end{equation}
Suppose that the extension \eqref{line bundle extension2} splits over $n$ fibres of $\pi$
(counting multiplicity). 
Referring to \eqref{long exact},
the first direct image sheaf $R^1\pi_\ast({K_1}^{-1} \otimes E)$
is then given by the extension 
\[ 0 \rightarrow \mathcal{S} \rightarrow R^1\pi_\ast({K_1}^{-1} \otimes E)
\rightarrow {H_+}^{-1} \rightarrow 0, \]
where $\mathcal{S}$ is a skyscraper sheaf supported on the $n$ points
(counting multiplicity) corresponding to these fibres.
By Grothendieck-
Riemann-Roch, 
\[ \deg(R^1\pi_\ast({K_1}^{-1} \otimes E)) = 
 -\frac{1}{2}c_1^2({K_1}^2 \otimes \delta^{-1}) = -\frac{1}{2} c_1^2(\pi^\ast (H_+))  = 0. \]
The degree of the line bundle $H_+$ is thus $n$; 
clearly, $H_+ = \mathcal{O}_B$ if $n=0$.
Note that by construction, $K_1 \otimes \delta^{-1} \otimes \pi^\ast(H_+) \otimes  
\mathcal{O}_X \left( \sum_{i=1}^r a_i T_i \right)$ is the destabilising line
bundle of $E$; however, 
$\pi_\ast(K_1 \otimes \delta^{-1} \otimes \pi^\ast(H_+) \otimes  
\mathcal{O}_X \left( \sum_{i=1}^r b_i T_i \right)) = \mathcal{O}_B$,
for all integers $0 \leq b_i \leq m_i - 1$.
Therefore, $a_i = m_i - 1$ for all $i = 1, \dots, r$, proving (ii).

Finally, let us assume that $\Sigma_1 \neq \Sigma_2$.
If the extension also splits over $m \leq 4\Delta(E)$ fibres of $\pi$
(counting multiplicity) corresponding to points in $\Sigma_1 \cap \Sigma_2$,
then the rank of $R^1\pi_\ast({K_1}^{-1} \otimes E)$
jumps at these $m$ points; in fact, the first direct image sheaf is given by the 
extension
\[ 0 \rightarrow \mathcal{O}_B \rightarrow R^1\pi_\ast({K_1}^{-1} \otimes E)
\rightarrow R^1\pi_\ast({K_1}^{-2} \otimes \delta) \rightarrow 0.\]
Dualising, we get $R^1\pi_\ast({K_1}^{-1} \otimes E)^\ast = H_-$, 
for $H_- \in \text{Pic}(B)$.
Let $n : = 4\Delta(E) - m$;
since the skyscraper sheaf $R^1\pi_\ast({K_1}^{-2} \otimes \delta)$ is supported
on $4\Delta(E)$ points (see the proof of Lemma \ref{not multiple section}), 
the line bundle $H_-$ has degree $-n$. 
Furthermore, by relative Serre duality,
$\pi_\ast(K_1 \otimes \delta^{-1} \otimes E) = R^1\pi_\ast({K_1}^{-1} \otimes E)^\ast = H_-$;
therefore, the second destabilising line bundle of $E$ is 
$K_2 = {K_1}^{-1} \otimes \delta \otimes H_-  \otimes \omega_{X/B}$ and we are done. 
\end{proof}

\begin{corollary}\label{filtrable bundle with trivial discriminant}
Every filtrable rank-$2$ bundle on $X$ that has no jumps and trivial discriminant is unstable.
\end{corollary}
\begin{proof}
Let $E$ be such a bundle;
its spectral cover of is then of the form $\left( \Sigma_1 + \Sigma_2 \right)$, 
for some sections $\Sigma_1,\Sigma_2 \subset J(X)$. 
Let $K_1$ be the destabilising bundle induced by $\Sigma_1$. 
Referring to Proposition \ref{max. destab. bundles}, 
if $\Sigma_1 = \Sigma_2$, then $\deg{K_1} \geq \deg{\delta}/2$ and $E$ is unstable.
If $\Sigma_1 \neq \Sigma_2$, then the extension \eqref{line bundle extension2} 
splits on every fibre of $\pi$
(because $\Delta(E) = 0$) and the second destabilising bundle
of $E$ is $K_2 = {K_1}^{-1} \otimes \delta \otimes \omega_{X/B}$; therefore,
$\deg{K_1} + \deg{K_2} = \deg{\delta} + \deg{\omega_{X/B}} \geq \deg{\delta}$ 
and at least one of the destabilising bundle
has degree greater or equal to $\deg{\delta}/2$.
\end{proof} 

Recall that, for surfaces $X$ with multiple fibres, 
the spectral cover of a vector bundle $E$ on $X$ was defined in section \ref{vector bundles} 
in terms of the vector
bundle $\psi^\ast E$ on an $m$-cyclic covering $\psi: X' \rightarrow X$, 
where $X'$ is an elliptic fibre bundle over a $m$-cyclic covering 
$B' \rightarrow B$.
Keeping this in mind, we now state the main result of the section.

\begin{theorem}\label{stable filtrable bundles}
Consider a filtrable rank-2 vector bundle $E$ on $X$ with determinant $\delta$
that has $k$ jumps of lengths $l_1, \dots , l_k$, respectively;
furthermore, suppose that $j$ of them occur over multiple fibres
$m_{i_1}T_{i_1}, \dots ,m_{i_j}T_{i_j}$, respectively, for some integer $0 \leq j \leq k$.
We set 
\[ \nu := \sum_{s=1}^j l_s / m_{i_s} + \sum_{t=j+1}^k l_t.\] 
Let $K$ be one of the destabilising bundles of $E$.
There is an extension 
\[ 0 \rightarrow \psi^\ast K \rightarrow \overline{\psi^\ast E} \rightarrow 
\psi^\ast(K \otimes \delta^{-1}) \rightarrow 0,\]
where $\overline{\psi^\ast E}$ denotes the vector bundle on $X'$ obtained by
performing successive elementary modifications to eliminate the jumps of $\psi^\ast E$.

(i) If $\Sigma_1 = \Sigma_2$ and the extension is trivial on every fibre of $\pi'$, 
then $E$ is stable if and only if $\nu > \deg{\omega_{X/B}}$.

(ii) Suppose that $\Sigma_1 = \Sigma_2$ and that the extension splits  
on only a finite number $mn$ of fibres, 
then $E$ is stable if and only if $\nu > n + \deg{\omega_{X/B}}$.

(iii) If $\Sigma_1 \neq \Sigma_2$ and the extension is non-trivial on
a finite number $mn$ of fibres of $\pi'$, then $E$ is stable if and only if
$\deg{K} \in \left( \deg{\delta}/2 - \nu - n + \deg{\omega_{X/B}},\deg{\delta}/2 \right)$.
\end{theorem}
\begin{proof}
Note that any elementary modification of $\psi^\ast E$ has the same destabilising 
bundles as $\psi^\ast E$ (which are the pullbacks to $X'$ of the destabilising bundles of $E$). 
Furthermore, the elementary modification $\overline{\psi^\ast E}$ 
has determinant 
$\psi^\ast \delta \otimes 
\mathcal{O}_{X'}(- D)$, where $D := \sum_{s=1}^j l_s T_{i_s} + (\sum_{t=j+1}^k l_t) T$.
Applying Proposition \ref{max. destab. bundles} to the bundle
$\overline{\psi^\ast E}$, we obtain the theorem.
\end{proof}

\section{Moduli spaces}

Let $X$ be a non-K\"{a}hler elliptic surface and consider 
a pair $(c_1,c_2)$ in $NS(X) \times \mathbb{Z}$.
For a fixed line bundle $\delta$ on $X$ with $c_1(\delta) = c_1$, let
$\mathcal{M}_{\delta,c_2}$ be the moduli space of stable holomorphic rank-2 vector bundles 
with invariants $\det(E) =  \delta$ and $c_2(E) = c_2$.
We define the following positive rational number:
\[ m(2,c_1) := -\frac{1}{4} \max 
\left\{ \sum_1^n \left( \frac{c_1}{2} - \mu_i \right)^2, 
\mu_1, \dots, \mu_r \in NS(X), \sum_1^n \mu_i = c_1 \right\} . \] 
\noindent
Note that, for any $c_1 \in NS(X)$, one can choose a line bundle $\delta$ on 
$X$ such that  
\begin{equation}\label{delta}
\mbox{$c_1(\delta) \in c_1 + 2NS(X)$ and 
$m(2,c_1) = -\dfrac{1}{2} \left( \dfrac{c_1(\delta )}{2} \right)^2$};
\end{equation} 
moreover, if there exist line bundles $a$ and $\delta'$ on $X$ such that 
$\delta = a^2\delta'$, then there is a natural isomorphism between the
moduli spaces $\mathcal{M}_{\delta,c_2}$ and $\mathcal{M}_{\delta',c_2}$,
defined by $E \mapsto a \otimes E$.
Therefore, if $\delta'$ is any other line bundle with Chern class in 
$c_1 + 2NS(X)$, it induces a moduli space that is isomorphic to 
$\mathcal{M}_{\delta,c_2}$. However, the advantage of using such a $\delta$ is that its 
Chern class has maximal self-intersection $-8m(2,c_1)$.
Hence, we restrict our study to moduli spaces $\mathcal{M}_{\delta,c_2}$
of stable bundles whose determinant $\delta$ satisfies \eqref{delta}.

\subsection{Existence and dimension}

A necessary condition for the existence of holomorphic rank-2 vector bundles is 
$\Delta(2,c_1,c_2) := 1/2 \left( c_2 - c_1^2/4 \right) \geq 0$ \cite{BaL,Br}. 
Also, a theorem of B\u{a}nic\u{a} - Le Potier's \cite{BaL} 
states that there exists a filtrable holomorphic rank-2 vector bundle with Chern classes $c_1$ and $c_2$ 
if and only if $\Delta(2,c_1,c_2) \geq m(2,c_1)$.
Given our choice of line bundle $\delta$, any element $E$ of $\mathcal{M}_{\delta,c_2}$ 
has discriminant 
\[ \Delta(E) = m(2,c_1) + \frac{1}{2}c_2  \geq 0.\]
Consequently, $c_2 \geq -2m(2,c_1)$; moreover, if $c_2 < 0$, then $E$ is unfiltrable.
However, if the vector bundle $E$ is an unfiltrable, then its spectral cover contains 
an irreducible bisection; it is then stable by Proposition \ref{spectral cover of stable bundles}.
Therefore, if a rank-2 vector bundle has second Chern class $-2m(2,c_1) \leq c_2 < 0$,
then it is stable.

Assume that the moduli space $\mathcal{M}_{\delta,c_2}$ is non-empty. Consider 
one of its elements $E$;
there is a natural splitting of the endomorphism bundle $\End(E) = \mathcal{O}_X \oplus 
\ad(E)$, where
$\ad(E)$ is the kernel of the trace map. By deformation theory,  
the moduli space has expected dimension $h^1(X;\ad(E)) - h^2(X;\ad(E))$ at $E$. 
Since the vector bundle $E$ is assumed to be stable, we have
$h^0(X;\ad(E)) = 0$ and the expected dimension of the moduli space is equal to 
$-\chi(E) = 8\Delta(2,c_1,c_2) - 3\chi(\mathcal{O}_X) = 8\Delta(2,c_1,c_2)$.

\subsection{Smoothness}
\label{dimension}

Let us first assume that $X$ is an elliptic fibre bundle
over a curve of genus less than $2$.
Recall that a vector bundle $E$ on a complex manifold $X$ is said to be {\em good} 
if and only if $h^2(X;\ad(E)) = 0$,
or equivalently, if $h^0(X;\ad(E) \otimes K_X) = 0$ (by Serre duality);
furthermore, the moduli space 
$\mathcal{M}_{\delta,c_2}$ is smooth at $E$ if and only if the vector bundle $E$ is good.
Given this, one easily proves the following:
\begin{proposition}
Let $X$ be a non-K\"{a}hler elliptic fibre bundle over
a curve $B$ of genus less than $2$, that is, $X$ is a Hopf surface or a primary Kodaira surface.
The moduli spaces $\mathcal{M}_{\delta,c_2}$ are then smooth of dimension $8\Delta(2,c_1,c_2)$.
\end{proposition}
\begin{proof}
It is sufficient to prove that every stable bundle on $X$ with Chern classes $c_1$ and $c_2$ is good.
In this case, the canonical bundle of the surface is $K_X = \pi^\ast(K_B)$. Since the genus of $B$ 
is $\leq 1$, the canonical bundle is given by $\mathcal{O}_X(-D)$, 
where $D$ is an effective divisor.
There is an inclusion $K_X = \mathcal{O}_X(-D) \subset \mathcal{O}_X$, 
which in turn induces an inclusion on the space of global sections 
$H^0(X;\ad{E} \otimes K_X) \subset H^0(X;\ad(E))$.
However, the stability of $E$ implies that $h^0(X;\ad(E)) = 0$ and we are done. 
\end{proof}

For an arbitrary non-K\"{a}hler elliptic surface 
$X \stackrel{\pi}{\rightarrow} B$,
we consider the elements of the moduli space that are regular,
that is, vector bundles that are regular on every fibre of $\pi$.
Note that for such a bundle $E$, the direct image sheaves 
$\pi_\ast(\End(E))$ and $R^1\pi_\ast(\End(E))$ are dual
locally free sheaves of rank two; therefore,
by Grothendieck-Riemann-Roch, we have
$c_1(\pi_\ast(\End(E))) = 2ch_2(E)$. 
Given the natural splitting $\pi_\ast(\End(E)) = \mathcal{O}_B \oplus \pi_\ast(\ad(E))$, 
we conclude that 
\[ \deg(\pi_\ast(\ad(E))) = 2ch_2(E) . \]
The Leray spectral sequence gives us 
$h^0(X;\ad(E) \otimes K_X) = h^0(B;\pi_\ast(\ad(E)) \otimes K_B)$; 
hence, if the degree of $\pi_\ast(\ad(E)) \otimes K_B$ is negative, 
we have $h^0(X;\ad(E) \otimes K_X) = 0$, 
leading us to the following: 
\begin{proposition}
Let $X$ be a non-K\"{a}hler elliptic surface over a base curve $B$ of genus $g$. 
Then, if $c_2 - c_1^2/2 > g - 1$, the moduli space $\mathcal{M}_{\delta,c_2}$
is smooth on the open dense subset of regular bundles.
\hfill \qedsymbol
\end{proposition}

\begin{remark}
We can also give a sufficient condition for smoothness of the moduli space at points that
do not correspond to regular bundles.
Consider a stable vector bundle $E$ that is not regular over the fibres of $\pi$
lying over the points $x_1, \dots, x_s$ in $B$. 
In this case, $\pi_\ast(\End(E))$ is again a rank-2 vector bundle, but $R^1\pi_\ast(\End(E))$ is 
the sum of a rank 2 vector bundle with a skyscraper sheaf supported on the points $x_1, \dots, x_s$, 
with multiplicities $\gamma_1, \dots, \gamma_s$, respectively. Let $\gamma  = \sum_i \gamma_i$.
Then, one easily verifies that a sufficient condition
for smoothness of the moduli space $\mathcal{M}_{\delta,c_2}$ at $E$ is given by
\[ c_2 - \frac{c_1^2}{2} > g-1 + \frac{\gamma}{4}.\]
Note that $\gamma$ depends not only on the spectral cover of $E$, but also on the 
geometry of its jumps.
\end{remark}

\subsection{The image of the graph map}
\label{image of graph map}
Fix any pair $(c_1,c_2) \in NS(X) \times \mathbb{Z}$ such that 
$\Delta(2,c_1,c_2) \geq 0$ and let $\delta$ be a line bundle on $X$ 
such that $m(2,c_1) = -\frac{1}{2} \left( c_1(\delta )/2 \right)^2$.
Referring to section \ref{rank-2 vector bundles},
this line bundle determines an involution $i_\delta$ of the Jacobian surface 
and there is an associated ruled surface $\mathbb{F}_\delta := J(X)/i_\delta$; 
the quotient map is denoted $\eta : J(X) \rightarrow \mathbb{F}_\delta$.
Furthermore, to any rank-2 vector bundle $E$ on $X$ 
with determinant $\delta$ and second Chern class $c_2$, 
there corresponds a graph in $\mathbb{F}_\delta$.
It was shown in \cite{Brinzanescu-Moraru1} that these graphs are 
elements of linear systems in $\mathbb{F}_\delta$ of the form $|\eta_\ast(B_0) + \mathfrak{b}f|$,
where $B_0$ is the zero section of $J(X)$, $\mathfrak{b}$ is the pullback to $X$ of a line bundle 
on $B$ of degree $c_2$, and $f$ is a fibre of the ruled surface.

Let $\mathbb{P}_{\delta,c_2}$ be the set of divisors in $\mathbb{F}_\delta$
of the form $\sum_{i=1}^k f_i + A$, where $A$ is a section 
and the $f_i$'s are fibres of the ruled surface,
that are numerically equivalent to $\eta_\ast(B_0) + c_2f$.  
We have a well-defined map
\[ G: \mathcal{M}_{\delta,c_2} \longrightarrow \mathbb{P}_{\delta,c_2} \]
that associates to each vector bundle its graph, called the
{\em graph map}.
Let us then describe the image of this map; 
we begin by noting that it 
is surjective on the open dense subset of graphs in $\mathbb{P}_{\delta,c_2}$ that
correspond to irreducible bisections in $J(X)$.
When considering the remaining graphs, we restrict ourselves, for simplicity, 
to the case where $X$ has no multiple fibres;
however, similar results holds if $X$ does have multiple fibres.
\begin{proposition}\label{image1}
Let $X \stackrel{\pi}{\rightarrow} B$ be a non-K\"{a}hler elliptic fibre bundle. 
Choose an element $c_1 \in NS(X)$ such that $m(2,c_1) = 0$; 
in this case, $\mathbb{F}_\delta = B \times \mathbb{P}^1$
and the elements of $\mathbb{P}_{\delta,c_2}$ are of the form
\[ \sum_{i=1}^k (\{ b_i \} \times \mathbb{P}^1) + Gr(F), \] 
where $b_1, \dots, b_k$ are points in $B$ and $Gr(F)$ is the graph of a rational map
$F: B \rightarrow \mathbb{P}^1$ of degree $c_2 - k$.
We then have the following.

(i) For $c_2 = 0$, the moduli space $\mathcal{M}_{\delta,0}$ is empty.

(ii) Let $\mathcal{S}$ be the set of points $\lambda_0$ in $T^\ast$ such that
the degree of any line bundle on $X$ corresponding to the section $B \times \{ \lambda_0 \}$
in $J(X)$ is congruent to $\deg{\delta}/2$ modulo $\mathbb{Z}$.
If $I$ is the projection of $\mathcal{S}$ onto $\mathbb{P}^1 = T^\ast/i_\delta$,
then we denote $B \times I$ the set of graphs
\[ \left\{ ( \{ b \} \times \mathbb{P}^1 ) + (B \times \{ \bar{\lambda} \}) \ \left| \ 
\mbox{$b \in B$ and $\bar{\lambda} \in I$} \right. \right\} . \]
For $c_2 = 1$, the image of the graph map $G$ is 
$\mathbb{P}_{\delta,1} \backslash \left( B \times I \right)$. 

(iii) For $c_2 \geq 2$, the graph map is surjective.
\end{proposition}
\begin{proof}
Consider a graph $\mathcal{G} = \sum_{i=1}^k (\{ b_i \} \times \mathbb{P}^1) + Gr(F)$, 
where $b_1, \dots, b_k$ are points in $B$ and $Gr(F)$ is the graph of a rational map
$F: B \rightarrow \mathbb{P}^1$ of degree $c_2 - k$;
we denote $\overline{C}$ the bisection of $J(X)$ determined by $Gr(F)$.
Referring to section \ref{vector bundles}, we can construct rank-2 vector bundles on $X$ 
with graph $\mathcal{G}$; therefore, we only have to determine whether or not 
at least one of them is stable.  
Let us fix a bundle $E$ with graph $\mathcal{G}$ and discuss its stability. 

Suppose that $c_2 = 0$; consequently, 
$\Delta(E) = 0$ and the map $F$ is constant.
Hence, the bisection $\overline{C}$ is reducible and $E$ is a filtrable bundle that
has no jumps and trivial discriminant; 
by Corollary \ref{filtrable bundle with trivial discriminant},
it is unstable, proving (i).

Now, assume that $c_2 \geq 1$. Recall that the vector bundle
$E$ may be unstable only if the bisection $\overline{C}$ is reducible; 
therefore, suppose that $\overline{C} = \Sigma_1 + \Sigma_2$ for some sections 
$\Sigma_1,\Sigma_2 \subset J(X)$.
If $\Sigma_1 = \Sigma_2$, then the bundle has at least one jump
(otherwise, $\Delta(E) = 0$ (see Lemma \ref{multiple section}), contradicting the fact that
$\Delta(E) = c_2 \geq 1$); in this case, $E$ is always stable by 
Theorem \ref{stable filtrable bundles}.
If $\Sigma_1 \neq \Sigma_2$, then a stable bundle $E$ can be constructed as follows.
Choose a line bundle $K$ corresponding to $\Sigma_1$;
after possibly tensoring $K$ by an element of $P_2$, one can assume that 
$\deg{K} \in \left( \deg{\delta}/2 - k - (c_2 - k), \deg{\delta}/2 \right)$,
unless $c_1 = 1$ and the degree of $K$ is congruent to $\deg{\delta}/2$ modulo $\mathbb{Z}$.
Then, consider a regular extension of $K^{-1} \otimes \delta(kT)$ by $K$
and perform $k$ elementary modifications (using a line bundle of degree $1$ on $T$) 
to introduce the jumps. Note that $K$ is one of the destabilising bundles of $E$;
referring to Theorem \ref{stable filtrable bundles}, $E$ is then stable.
Finally, if $c_1 = 1$ and $\Sigma_1 \neq \Sigma_2$, then a bundle with graph $\mathcal{G}$ is stable
if and only if the degrees of its destabilising bundles are in the interval
$\left( \deg{\delta}/2 - 1, \deg{\delta}/2 \right)$;
if all bundles corresponding to $\Sigma_1$ have degree congruent to 
$\deg{\delta}/2$ modulo $\mathbb{Z}$, then this is never possible.
\end{proof}

\begin{proposition}\label{image2}
Let $X \stackrel{\pi}{\rightarrow} B$ be a non-K\"{a}hler elliptic fibre bundle. 
Choose an element $c_1 \in NS(X)$ such that $m(2,c_1) > 0$, so that we may have $c_2 < 0$.
Then, the graph map is surjective whenever the moduli spaces are non-empty,
except in the following case.
Suppose that $c_2 = 0$ and $m(2,c_1) = 1/4$.
Furthermore, let $J$ be the set of sections $A$ in $\mathbb{P}_{\delta,0}$ such that 
$\eta^\ast A = \Sigma_1 + \Sigma_2$ is a reducible bisection of $J(X)$ and the degree
of any line bundle on $X$ associated to $\Sigma_1$ is 
congruent to $\deg{\delta}/2$ modulo $\mathbb{Z}$.
In this case, the image of the 
graph map is $\mathbb{P}_{\delta,0} \backslash J$.
\end{proposition}
\begin{proof}
Consider a graph $\mathcal{G} = \sum_{i=1}^k f_i + A$ and set $\overline{C} = \eta^\ast A$. 
We know that there exist bundles corresponding to this graph; 
let us then discuss the stability of a bundle $E$ that has graph $\mathcal{G}$.
If $c_2 < 0$, then $\Delta(E) < m(2,c_1)$: the bisection $\overline{C}$ is irreducible and 
the bundle is stable.
If $c_2 = 0$, then $\Delta(E) = m(2,c_1) > 0$. 
There are now two possibilities. The first is $k \neq 0$,
implying that $A^2 < 4m(2,c_1)$; therefore, the bisection $\overline{C} = \eta^\ast A$ is 
irreducible and the bundle is stable. The second is $k=0$ 
and the bisection is reducible;
suppose that $\overline{C} = \Sigma_1 + \Sigma_2$, for some sections $\Sigma_1, \Sigma_2 \subset J(X)$.
Note that $\Sigma_1 \neq \Sigma_2$; otherwise, $k=0$ would imply that $\Delta(E) = 0$, 
which is a contradiction.
The vector bundle $E$ is then an extension of $K^{-1} \otimes \delta$ by $K$,
where $K$ is the destabilising bundle of $E$ corresponding to $\Sigma_1$, that
can be assumed to be regular on every fibre of $\pi$. Hence, $E$ is stable if and only if
$\deg{K} \in \left( \deg{\delta}/2 - 4m(2,c_1), \deg{\delta}/2 \right)$,
as stated in Theorem \ref{stable filtrable bundles}.
Clearly, if $m(2,c_1) = 1/4$ and the degree of every line bundle corresponding to
$\Sigma_1$ is congruent to $\deg{\delta}/2$ modulo $\mathbb{Z}$, then $E$ is unstable.
Finally, by arguments similar to those used to prove Proposition \ref{image1},
the graph map is surjective whenever $c_2 \geq 1$.
\end{proof}

\subsection{Fibre of the graph map}
\label{general jumps}

If we consider graphs without vertical components, the description
of most fibres of the graph map is then straightforward.

\begin{proposition}\label{generic fibre}
Let $X$ be a non-K\"{a}hler elliptic surface over a curve $B$ of genus $g$.
Fix a pair $(c_1,c_2)$ in $NS(X) \times \mathbb{Z}$ and let $\delta$ be a line bundle on $X$,
with $c_1(\delta) = c_1$, such that $\mathcal{M}_{\delta,c_2}$ is non-empty. 
Consider an element $A$ of $\mathbb{P}_{\delta,c_2}$ that does not contain vertical components
and let $\overline{C} = \eta^\ast A$ be the corresponding bisection in $J(X)$.

(i) Suppose that $\overline{C}$ is a smooth bisection of $J(X)$. 
The fibre of the graph map $G$
at $A$ is then isomorphic to a finite number of copies of the Prym variety $Prym(\overline{C}/B)$.
Note that if the surface $X$ does not have multiple fibres, then 
$G^{-1}(A)$ consists of a single copy of $Prym(\overline{C}/B)$
(see Theorem \ref{classification of regular bundles}).

(ii) If the bisection $\overline{C} = \Sigma_1 + \Sigma_2$ is reducible,
then the components of $G^{-1}(A)$ are 
parametrised by the set of line bundles on $X$ associated to $\Sigma_1$
that satisfy the conditions of
Theorem \ref{stable filtrable bundles}.
In particular, the component given by the line bundle $K$ consists of 
extensions of $K^{-1} \otimes \delta$ by
$K$ that are regular on at least one fibre $\pi^{-1}(b)$ where
$\Sigma_{1,b} = \Sigma_{2,b}$, if $\deg{K}$ is not congruent to $\deg{\delta}/2$ 
modulo $\mathbb{Z}$,
or that are regular on at least two such fibres, otherwise. 
\hfill \qedsymbol
\end{proposition}

For graphs with vertical components, the fibre of the graph map 
can be described by examining how jumps can be added to vector bundles,
that is, by classifying elementary modifications.
This is done in detail in \cite{Moraru} for vector bundles on Hopf surfaces.
For the sake of completion, we briefly state how this translates to bundles
on an arbitrary non-K\"{a}hler elliptic surface $X$.

Let $E$ be a stable rank-2 vector bundle on $X$ with $\det{E} = \delta$, $c_2(E) = c_2$, 
and a jump of length $l$ over the smooth fibre $T = \pi^{-1}(x_0)$.
This jump can be removed by performing $l$ successive allowable
elementary modifications, thus obtaining a bundle with determinant $\delta(-lT)$;
note that this procedure is canonical.
But, adding a jump to $E$ implies several choices: a jumping sequence
$\{ h_0, \dots, h_{l-1} \}$, a line bundle $N$ on $T$ for each distinct integer of the
jumping sequence, and surjections to $N$ that preserve stability.
These choices are parametrised by a fibration that we now describe. 
Let $\mathcal{G}$ be a graph that contains a vertical component over $x_0$ of
multiplicity $\mu$ and $\{ h_0,\dots,h_{l-1} \}$ be a jumping sequence
such that $\sum_{i=0}^{l-1} h_i = \mu$. We set
\[ \mathcal{E}^j \mathcal{J}^{c_2,l}_{\mathcal{G},\{ h_0,\dots,h_{l-1} \}} =
\left\{ E \in \mathcal{M}_{\delta(jT),c_2}  \ \left|
\begin{array}{c}
\mbox{$G(E) = \mathcal{G}$ and $E$ has a jump}\\
\mbox{of length $l$ at $x_0$ with jumping }\\
\mbox{sequence $\{ h_0,\dots,h_{l-1} \}$}
\end{array} \right. \right\}. \] 
Associating to a bundle $E$ its allowable elementary modification
$\bar{E}$ therefore defines a natural map
\[ \begin{array}{rcl}
\Psi : \mathcal{E}^{j+1} \mathcal{J}^{c_2,l+1}_{\mathcal{G},\{ h_0,h_1,\dots,h_{l-1} \}} & 
\longrightarrow &
\mathcal{E}^j \mathcal{J}^{c_2-h_0,l}_{G(\bar{E}),\{h_1,\dots,h_{l-1} \}} \\
E & \longmapsto & \bar{E}. \end{array} \]
\begin{proposition}
The fibre of the natural projection $\Psi$ at $W$ is given by:

(i) $\text{Aut}_{\text{SL}(2,\mathbb{C})}(W|_T)$, \mbox{if $c_2>h_0$ and $h_0 = h_1$},

(ii) $\displaystyle \text{Pic}^{-h_0}(T) \times \text{Aut}_{\text{SL}(2,\mathbb{C})}(W|_T)$,
\mbox{if $c_2 > h_0 > h_1$ and $l>0$},

(iii) $\displaystyle \text{Pic}^{-c_2}(T)$, \mbox{if $c_2=h_0$ or $l=0$}. 
\hfill \qedsymbol
\end{proposition}

\subsection{Integrable systems}

A Poisson structure on a surface $X$ is given by a global 
section of its anticanonical bundle $K_X^{-1}$ \cite{Bottacin-Poisson}.
Suppose that $X \stackrel{\pi}{\rightarrow} B$ is a non-K\"{a}hler elliptic surface
that may have multiple fibres $T_1, \dots , T_r$ of multiplicities $m_1, \dots , m_r$,
respectively. 
The anticanonical bundle of $X$ is then
$\pi^\ast K_B^{-1} \otimes \left( \bigotimes_{i=1}^r \mathcal{O}_X(1 - m_i) \right)$,
implying that $X$ admits a Poisson structure if and only if the genus of the base curve is
$\leq 1$ and if the surface does not have multiple fibres.
From now on, we suppose that $X$ is a non-K\"{a}hler elliptic surface without multiple fibres
over a curve $B$ of genus $g = 0$ or $1$, that is, a Hopf surface or a primary Kodaira surface.
Let us fix a Poisson structure $s \in H^0(X,K_X^{-1})$ on $X$. 
A Poisson structure $\theta = \theta_s \in 
H^0(\mathcal{M},\otimes^2T\mathcal{M})$ on the moduli space 
$\mathcal{M} := \mathcal{M}_{c_2,\delta}$ is then defined as follows: for any bundle
$E \in \mathcal{M}$, 
$\theta(E) : T^\ast_E\mathcal{M} \times T^\ast_E\mathcal{M}
\longrightarrow \mathbb{C}$ 
is the composition
\[ \begin{array}{lr}
\theta(E): H^1(X,\ad(E) \otimes K_{X}) \times
H^1(X,\ad(E) \otimes K_{X}) 
\stackrel{\circ}{\longrightarrow} \\
\hspace{0.9in} H^2(X,\End(E) \otimes K^2_{X}) 
\stackrel{s}{\longrightarrow}
H^2(X,\End(E) \otimes K_{X}) 
\stackrel{\text{Tr}}{\longrightarrow} \mathbb{C},
\end{array} \]
where the first map is the cup-product of two cohomology classes, the second is 
multiplication by $s$, and the third is the trace map. 

If the base curve $B$ is elliptic, the canonical bundle of $X$ is trivial and the Poisson 
structure $s$ is non-degenerate; in this case, $\theta$ has maximal rank everywhere,
that is, $\theta$ is symplectic.
If the base curve is instead rational, 
the Poisson structure $s$ is now degenerate; we denote its divisor $D := (s)$.  
Then, at any point $E \in \mathcal{M}$, 
\[ \rk \theta(E) = 4 \dim_\mathbb{C} \mathcal{M}  - \dim H^0(D,\ad(E|_D)). \]
Suppose that the locally free sheaf $\mathcal{O}_B(2)$  on $B \cong \mathbb{P}^1$ is given
by the divisor $x_1 + x_2$, for some points $x_1, x_2 \in B$;
then, $D = T_1 + T_2$, where $T_i=\pi^{-1}(x_i)$ for $i=1,2$.
We now see that the rank of the Poisson structure is generically $4\dim_\mathbb{C} \mathcal{M} -2$
and ``drops'' at the points of $\mathcal{M}$ corresponding to bundles that 
are not regular over the fibres $T_1$ and $T_2$ (for details, see \cite{Moraru}).
 
Referring to sections \ref{dimension} and \ref{general jumps}, 
the moduli space $\mathcal{M}$ has dimension $8\Delta(2,c_1,c_2)$
and the generic fibre of the graph map 
$G: \mathcal{M} \longrightarrow \mathbb{P}_{\delta,c_2}$
is a Prym variety of dimension $4\Delta(2,c_1,c_2) + g - 1$
(see Proposition \ref{generic fibre}).
Also, one can show as in \cite{Moraru} that
the component functions $H_1, \dots, H_N$ of the graph map are in involution
with respect to the Poisson structure, 
that is, $\{ H_i,H_j \} = 0$ for all $i,j$.
Consequently, the graph map $G$ is an
algebraically completely integrable Hamiltonian system.

\end{document}